\documentclass{article}%
\usepackage{amsmath, amsthm}
\usepackage{amsfonts}
\usepackage{amssymb}
\usepackage{enumerate}
\usepackage{graphicx}
\usepackage[table]{xcolor}
\usepackage[all,arc,curve,color,frame]{xy} % til at tegne grafer og diagrammer

\usepackage{amsmath}
\usepackage{amsfonts}
\usepackage{amssymb}
\usepackage{graphicx}%
\usepackage{color}
\usepackage{amsmath}
\usepackage{amsfonts}
\usepackage{amssymb}
\usepackage{graphicx}%
\usepackage{xypic}
\usepackage{tikz-cd}
\usepackage{color}
\usepackage{mathbbol}
\providecommand{\U}[1]{\protect\rule{.1in}{.1in}}
\newtheorem{theorem}{Theorem}

\newtheorem{corollary}[theorem]{Corollary}

\newtheorem{lemma}[theorem]{Lemma}

\newtheorem{proposition}[theorem]{Proposition}
\newtheorem{remark}[theorem]{Remark}

\tikzset{
  symbol/.style={
    draw=none,
    every to/.append style={
      edge node={node [sloped, allow upside down, auto=false]{$#1$}}}
  }
}

\title{Sections and cones}

\author{Tatiana Shulman}

\begin{document}

%\author{.......}
%\address{Department of Mathematical Sciences, UNiversity of Gothenburg, Chalmers tvärgata 3, 412 96 Gothenburg, Sweden}
%\email{tatshu@chalmers.se}

%\author{}
%\address{}
%\email{}

\maketitle

\begin{abstract} By Bartle-Graves theorem every surjective map between C*-algebras has a continuous section, and Loring proved that there exists a continuous section of norm arbitrary close to 1. Here we prove that there exists a continuous section of norm exactly 1. This result is used in the second part of this paper which is devoted to properties of cone C*-algebras.

  We prove that any $\ast$-homomorphism from the cone over a separable C*-algebra to a quotient C*-algebra always lifts to a contractive asymptotic homomorphism.
 As an application we give a short proof and a strengthening of the result of Forough-Gardella-Thomsen that states that any cpc (order zero) map  has an asymptotically cpc  (order zero, respectively) lift.

 As another application we prove that all hyperlinear traces on cones and suspensions are MF.
We also give unified proofs of Voiculescu's result that cones are quasidiagonal and Brown-Carrion-White's result  that all amenable traces on cones are quasidiagonal.
\end{abstract}

%\subjclass[2010]{Primary 46L05; Secondary  20E26,  46L55}

%\keywords{}

\tableofcontents

\bigskip

\bigskip

\section{Introduction}

Lifting $\ast$-homomorphisms is crucial in many aspects of  C*-theory.  It is desirable but rather rare that a $\ast$-homomorphism lifts to a $\ast$-homomorphism. However in many situations it suffices to lift $\ast$-homomorphisms to maps from broader classes. In the series of works on extensions of C*-algebras Manuilov and Thomsen consider liftings to asymptotic homomorphisms, that is maps that are only asymptotically multiplicative, asymptotically linear etc. (see e.g. \cite{MT15} and references therein).  Still, not all C*-algebras allow such lifting \cite{MT2}.  On the other hand, if one drops all the "algebraic" requirements and asks only about a continuous lift, then it is always possible to find such. This is stated by  Bartle-Graves theorem and is exceptionally useful in C*-theory. A section whose existence is guaranteed by Bartle-Graves theorem is not linear. However  one often needs linear continuous or, even better,  completely positive lifts of $\ast$-homomorphisms. As is well known, this is possible for  many C*-algebras (e.g. for all nuclear C*-algebras by Choi-Effros theorem) but not for all. It is therefore surprising that any completely positive contractive map admits an asymptotically completely positive asymptotically contractive lift, as was discovered by Forough-Gardella-Thomsen \cite{FGT}. They also obtained an equivariant version and  an order zero version of their lifting result.

\medskip

The first part of this paper deals with continuous sections. While Bartle-Graves theorem does not say anything about the norm of a continuous section,  Loring proved  that for a surjective $\ast$-homomorphism between C*-algebras there exists a continuous section of norm arbitrary close to one (\cite{LoringSection}, \cite[Th. 3.2]{Manuilov}). Here we prove that there exists a continuous section of norm exactly one. Moreover the section can be made positive.

\medskip

{\bf Theorem.} {\it (Theorem \ref{section}) Any surjective $\ast$-homomorphism $\pi: B \to D$ has a continuous section $s: D \to B$ which is
\begin{enumerate}

\item  contractive,

\item  unital, if $\pi$ if unital,

\item selfadjoint (that is $s(a^*) = s(a)^*$ for any $a\in D$),

\item positive (that is $s(a)\ge 0$ for any $a\ge 0$).
\end{enumerate}

\noindent If $G$ is  a compact group, $B, D$ are $G$-algebras and $\pi$ is equivariant, then $s$  can be chosen equivariant.
}

\medskip

\noindent We also prove that given a continuous  (contractive, etc.) section on a closed subset of a C*-algebra  one always can extend it to a continuous  (contractive, etc.) section on the whole C*-algebra ({\bf Theorem \ref{ExtendingSection}}).

These results are used in the second part of the paper.

The second part of this work is about properties of cone C*-algebras. In \cite{LoringShulman} Loring and the author  found that cones are generated by homogeneous relations. It was one of two ingredients that allowed us to prove that cones over separable C*-algebras can be written as inductive limits of projective C*-algebras \cite{LoringShulman}. In the second part of the present work we demonstrate that  the fact  that cones are generated by homogeneous relations alone leads to obtaining various other (known and unknown) statements about cones.  The main results of the second part of this paper are the following two theorems.

\medskip

{\bf Theorem.} (Theorem \ref{cones}) {\it Let $A$ be a  separable C*-algebra. Then any $\ast$-homomorphism from $CA$ to any quotient  lifts to a contractive positive asymptotic homomorphism.}

\medskip

The second result is a completely positive version of the theorem above. Its proof required developing some additional technique.

\medskip

{\bf Theorem.} {\it (Theorem \ref{ccp})   Let $A$ be separable and let $f: CA \to B/I$ be a $\ast$-homomorphism that has a completely positive contractive (cpc, for short) lift $\phi: CA \to B$.
 Then $f$ lifts to a cpc  asymptotic homomorphism.}

  \medskip

These results have several applications. First of all, we give a short proof of the aforementioned result of Forough - Gardella - Thomsen \cite{FGT}. Our proof shows additionally that an asymptotically contractive asymptotically cp (asymptotically order zero, respectively) lift of a cpc (order zero, respectively) map can be chosen to be the composition of an order zero map and a corner of an asymptotic homomorphism (an asymptotic homomorphism, respectively). Moreover the lift can be made not only asymptotically contractive but exactly contractive and exactly positive. Equivariant version of it is also given. ({\bf Corollary \ref{ShortProof}, Corollary \ref{corollaryOrderZero}, Proposition \ref{equivariant}}).

The next applications deal with traces on C*-algebras. Approximation properties of traces -- quasidiagonality, amenability, hyperlinearity and the property of being an MF-trace -- nowadays play crucial role in the classification program for simple nuclear C*-algebras (\cite{TWW},  \cite{Gabe}) and at the study of group C*-algebras (see e.g. \cite{RainoneSch} and \cite{SchafhauserAmalgamated}).
Whether all amenable traces are quasidiagonal and all hyperlinear traces are MF is an open problem.

As another application of the theorems above  we obtain a unified proof of two well-known results about cones: cones are quasidiagonal (Voiculescu \cite{Voiculescu}) and all amenable traces on cones are quasidiagonal (Brown-Carrion-White \cite{BCW}), see {\bf Theorems \ref{QDtraceCone} } and {\bf \ref{Voiculescu}}. We also prove the following new fact.

\medskip

{\bf Theorem.} {\it (Theorem \ref{MFtraceCone}) All hyperlinear traces on cones are MF. }

\medskip

\noindent The same statement is proved for suspensions ({\bf Th. \ref{suspension}}).

\medskip

As byproduct of this work we obtain two results that can be of independent interest:

1) By Winter-Zacharias duality, order zero maps $\rho: A \to B$ are in one-to-one correspondence with $\ast$-homomorphisms  $\rho_{\phi}: CA \to B$ \cite{WZ}. As byproduct of techniques developed in this paper, we obtain the following

\medskip

{\bf Corollary.} {\it  (Corollary \ref{OrderZeroVsHomomorphism}) Let $A$  be a separable C*-algebra and $\phi:A \to B/I$  an order zero map. If $\phi $ lifts to a cp map, then $\rho_{\phi}$ lifts to a cp map.}

\medskip

\noindent As a particular case it covers \cite[Lemma 4.1]{BCW} with a short proof.  \footnote{\cite[Lemma 4.1]{BCW} is exactly what made the proof of quasidiagonality of amenable traces on cones long in \cite{BCW}.}

\medskip

2) The Lifting Property (LP) and the Local Lifting Property (LLP) are introduced by Kirchberg  in his seminal paper \cite{Kirchberg}.
 The (L)LP deals with lifting of completely positive maps and  has deep connections with tensor products and the weak expectation property of C*-algebras  and is one of key ingredients in the proving the equivalence of several fundamental conjectures.

 \noindent While the LLP is known to be preserved under extensions (\cite[Cor. 2.6.]{Kirchberg}), it appears  to be unknown for the LP. We obtain it here as a byproduct of this work.

{\bf Corollary.} {\it (Corollary \ref{LP}) If $0\to J\to E\to A\to 0$ is a short exact sequence,  $E$ is separable, and $J$ and $A$ have the Lifting Property (LP), then $E$ has the LP.}

\medskip

\noindent {\bf Acknowledgments. } The author is grateful to Dominic Enders for  very useful comments and suggestions on a draft of this paper.
  The author is partially supported by a grant from the Swedish Research Council no. 2025-05359. The updated version of this paper is done while the author was in residence at Institut Mittag-Leffler in Djursholm during the semester ”Operator algebras and Quantum information theory” Spring, 2026 supported by the Swedish Research Council under grant no. 2021-06594.

\section{Preliminaries}

{\bf Asymptotic homomorphisms. }

\medskip

{\bf Definition} (\cite{ConnesHigson}). An {\it asymptotic homomorphism} from $A$ to $B$ is a family of maps $(f_{\lambda})_{\lambda\in [0,\infty)}: A \to B$ satisfying the following properties:

\medskip

(i) for any $a\in A$, the mapping $[0, \infty) \to B$ defined by the rule $\lambda \to f_{\lambda}(a)$  is continuous;

(ii) for any $a, b\in A$ and $\mu_1, \mu_2\in \mathbb C$, we have

\begin{itemize}
\item  $\lim_{\lambda\to \infty} \|f_{\lambda}(a^*) - f_{\lambda}(a)^*\| = 0;$
\item  $\lim_{\lambda\to \infty} \|f_{\lambda}(\mu_1a + \mu_2b) - \mu_1f_{\lambda}(a) - \mu_2f_{\lambda}(b)\| = 0;$
\item $\lim{\lambda\to \infty}  \|f_{\lambda}(ab) - f_{\lambda}(a)f_{\lambda}(b)\| = 0.$
\end{itemize}

We will call $(f_{\lambda})_{\lambda\in \Lambda}: A \to B$, where $\Lambda$ is a directed set,   a {\it discrete asymptotic homomorphism} if the condition (ii) above is satisfied and for each $a\in A$ one has $\sup_{\lambda} \|f_{\lambda}(a)\|< \infty$. (For usual asymptotic homomorphisms the last condition holds automatically, see \cite{ConnesHigson} or \cite{Dadarlat2}). In this paper discrete asymptotic homomorphisms mostly will be indexed by $\Lambda = \mathbb N$.

Two (discrete) asymptotic homomorphisms $(f_{\lambda})_{\lambda\in [0, \infty)}, (g_{\lambda})_{\lambda\in \Lambda}: A \to B$ are {\it equivalent} if for any $a\in A$, we have $\lim_{\lambda\to \infty} \|f_{\lambda}(a)- g_{\lambda}(a)\| = 0.$

Two (discrete) asymptotic homomorphisms $(f_{\lambda})_{\lambda\in \Lambda}, (g_{\lambda})_{\lambda\in \Lambda}: A \to B$ are {\it homotopy equivalent} if there exists an asymptotic homomorphism $(\Phi_{\lambda})_{\lambda\in \Lambda}: A \to B\otimes C[0, 1]$ such that $ev_0\circ \Phi_{\lambda} = f_{\lambda}, \; ev_1\circ \Phi_{\lambda} = g_{\lambda},$ $\lambda\in \Lambda$.

A (discrete) asymptotic homomorphism $(f_{\lambda})_{\lambda\in \Lambda}$ is {\it equicontinuous} if for any $a_0\in A$ and $\epsilon >0$ there is a $\delta>0$ such that $\|f_{\lambda}(a) - f_{\lambda}(a_0)\|<\epsilon$ for every $\lambda\in \Lambda$ whenever $\|a-a_0\|\le\delta$.

\medskip

Following Manuilov and Thomsen we will say that a $\ast$-homomorphism $f: A\to B/I$ lifts to an asymptotic homomorphism  $\phi_t$, $t\in [1, \infty)$,  if $q\circ \phi_t = f$, for each $t$.

\medskip

%Asymptotic lifting. A $\ast$-homomorphism $f: A\to B/I$ asymptotically lifts to an asymptotic homomorphism  $\phi_t$ if $\lim_{t\to \infty} q\circ \phi_t(a) = f(a)$, for any $a\in A$. Carrion and Schafhauser considered such liftings in \cite{CSch}.

%\begin{lemma}\label{AsMultiplicativeFamily} Suppose $f_{\lambda}: C^*\langle \mathbf x\;|\;\mathbf R\rangle \to B$, $\lambda\in \Lambda$, is a family of maps such that $$\lim_{\lambda\to\infty}\|\mathbf R(f_{\lambda}(\mathbf x))\|= 0.$$ Then the family $f_{\lambda}, \lambda\in \Lambda$, is asymptotically multiplicative, asymptotically linear and asymptotically self-adjoint.\end{lemma}
%\begin{proof} Let $\pi: \prod_{\lambda} B \to \prod_{\lambda} B/\bigoplus_{\lambda} B$ be the canonical surjection. Then $\pi\circ(f_{\lambda})_{\lambda\in \Lambda}$ is a $\ast$-homomorphism. Therefore, for any $a, b\in C^*\langle \mathbf x\;|\;\mathbf R\rangle$,$$\pi\left(\left(f_{\lambda}(ab) - f_{\lambda}(a)f_{\lambda}(b)\right)_{\lambda\in \Lambda}\right)=0$$ that means that $$\lim_{\lambda\to\infty}\|f_{\lambda}(ab) - f_{\lambda}(a)f_{\lambda}(b)\|= 0.$$ Similarly one shows asymptotic linearity and self-adjointness.\end{proof}

We will call a (not necessarily linear) map {\it positive} if  it sends positive elements to positive elements.

%$\prod M_n/\oplus_{2, \omega} M_n$. We denote the trace .... on it by $tr$.

\bigskip

{\bf Cones.}

\medskip

We will say that a non-commutative $\ast$-polynomial $p(x_1, \ldots, x_n)$ is {\it $d$-homogeneous} if $p(tx_1, \ldots, tx_n) = t^d p(x, y)$
for all real scalars $t$.  We call $d$  the {\it degree of homogeneity} of $p$.

For an NC $\ast$-polynomial $p(x_1, \ldots, x_n)$,   its {\it homogenization} $\tilde p(h, x_1, \ldots, x_n)$ is the homogeneous NC $\ast$-polynomial derived from $p$ by padding monomials on the left with
various powers of $h$. For example, if $$p(x_1, x_2, x_3) = x_1^4x_3 - x_2^*x_1 + x_1^*,$$ then
$$\tilde p (h, x_1, x_2, x_3) = x_1^4x_3 - h^3x_2^*x_1 + h^4x_1^*.$$

Let $A$ be a separable unital C*-algebra given by presentation
\begin{equation}\label{PresentationATAMS}A = \left\langle x_1, x_2, \ldots\;|\; -c_i\le x_i\le c_i,\; R_k(x_1, x_2, \ldots) =0, i, k\in \mathbb N\right\rangle\end{equation}.  Here and everywhere throughout this paper each
NC $\ast$-polynomial depends only of finitely many variables. By \cite[Lemma 7.3]{LoringShulman} any separable C*-algebra is of the form (\ref{PresentationATAMS}).  In \cite[Lemma 7.1]{LoringShulman} the following presentation  for the cone $CA$  was found
%\footnote{ In \cite[Lemma 7.1]{LoringShulman} there was assumed that the relations of $A$ do  not  have free term, but in fact this requirement can be omitted. }
:

\begin{equation}\label{ConePresentationTAMS}  CA = \left\langle  h, x_1, x_2, \ldots\;|\; 0\le h\le 1, -c_ih\le x_i\le c_ih, [h, x_i] = 0, \tilde R_k(h, x_1, x_2, \ldots) = 0, i, k \in \mathbb N \right\rangle \end{equation}

Thus the cone has a presentation where every relation is homogeneous, the fact that will be crucial for this paper.

We note that if $A$ is non-unital, then, as the proof of \cite[Lemma 7.1]{LoringShulman} shows, all the elements $x_i$ and $x_ih^k$, $k\in \mathbb N$,  belong to $CA$ and generate it.

\bigskip

{\bf Traces}

\medskip

Recall that a trace $\tau$ on a C*-algebra $A$ is {\it amenable} if there is a sequence of   ccp maps $\phi_n: A \to M_{k_n}$, $n\in \mathbb N$, such that
$$ \lim_{n\to \infty} \|\phi_n(ab)- \phi_n(a)\phi_n(b)\|_2=0 \;\; \text{and} \;\; \lim_{n\to \infty}|\tau(a)- tr \phi(a)|=0,$$ for any $a, b\in A$.

A trace $\tau$ is {\it quasidiagonal} if there is a sequence of   ccp maps $\phi_n: A \to M_{k_n}$, $n\in \mathbb N$, such that
$$ \lim_{n\to \infty} \|\phi_n(ab)- \phi_n(a)\phi_n(b)\|=0 \;\; \text{and} \;\; \lim_{n\to \infty}|\tau(a)- tr \phi(a)|=0,$$  for any $a, b\in A$.

If in the the definition of an amenable trace we drop the requirement that  $\phi_n$'s are ccp, then we obtain the notion of a hyperlinear trace. Precisely, a trace $\tau$ is {\it hyperlinear} if there is a sequence of  maps $\phi_n: A \to M_{k_n}$, $n\in \mathbb N$, such that
\begin{multline*}  \lim_{n\to \infty} \|\phi_n(ab)- \phi_n(a)\phi_n(b)\|_2=0 \;\; \lim_{n\to \infty}\|\phi_n(\lambda a+\mu b)- \lambda\phi_n(a)-\mu \phi_n(b)\|_2 =0, \\ \lim_{n\to \infty} \|\phi_n(a^*) - \phi_n(a)^*\|_2=0,\;\; \sup_{n\in \mathbb N} \|\phi_n(a)\|<\infty, \;\text{and} \;\; \lim_{n\to \infty}|\tau(a)- tr \phi(a)|=0,\end{multline*} for any $a, b\in A$, $\lambda, \mu\in \mathbb C$.

 \begin{remark}\label{hypTraceReformulation} Equivalently, one can say that $\tau$ is {\it hyperlinear}  if $\tau = tr\circ f$, for some $\ast$-homomorphism $f: A \to \prod M_{k_n}/\oplus_2 M_{k_n}$, for some $k_n$'s,  where $\oplus_2 M_{k_n}$ is the ideal of all sequences that converge to zero in the 2-norm and $tr$ is a  trace on $\prod M_{k_n}/\oplus_2 M_{k_n}$ defined by the formula
 $tr([(T_n)_{n\in \mathbb N}]) = \lim_{n\to \omega} tr T_n$, where $\omega$ is some non-trivial ultrafilter on $\mathbb N$. One also can consider the ideal $\oplus_2 M_{k_n} = \{(T_n)_{n\in \mathbb N}\;|\; \lim_{n\to \omega} \|T_n\|_2 = 0\}$ of $\prod M_{k_n}$. The C*-algebra $\prod M_{k_n}/\oplus_{2, \omega} M_{k_n}$, called  the tracial ultraproduct of matrix algebras,  is a von Neumman algebra with the faithful trace $tr([(T_n)_{n\in \mathbb N}]) = \lim_{n\to \omega} tr T_n$. It is easy to show that a trace $\tau$ is hyperlinear if and only if $\tau = tr\circ f$, for some $\ast$-homomorphism $f: A \to \prod M_{k_n}/\oplus_{2, \omega} M_{k_n}$ and some $k_n$'s.

\noindent One can also replace  $\prod M_n/\oplus_{2, \omega} M_{n}$ by $\mathcal R^{\omega}$.
\end{remark}

If in the definition of a quasidiagonal trace one drops the requirement that $\phi_n$'s are ccp, one obtains the notion of an MF-trace. Precisely,
a trace $\tau$ is {\it MF} if there is a sequence of  maps $\phi_n: A \to M_{k_n}$, $n\in \mathbb N$, such that
\begin{multline*}  \lim_{n\to \infty} \|\phi_n(ab)- \phi_n(a)\phi_n(b)\|=0 \;\; \lim_{n\to \infty}\|\phi_n(\lambda a+\mu b)- \lambda\phi_n(a)-\mu \phi_n(b)\| =0, \\ \lim_{n\to \infty} \|\phi_n(a^*) - \phi_n(a)^*\|=0,\;\; \sup_{n\in \mathbb N} \|\phi_n(a)\|<\infty, \; \text{and} \;\; \lim_{n\to \infty}|\tau(a)- tr \phi(a)|=0,\end{multline*} for any $a, b\in A$, $\lambda, \mu\in \mathbb C$.

\begin{remark} In the definition of an MF-trace  one can additionally require the maps $\phi_n$ to be $\ast$-linear   \cite[Prop. 2.2]{RainoneSch}. The same argument  works also for  hyperlinear traces.
%Indeed the homomorphism $f: A \to \prod M_n/\oplus M_n$ can be lifted to a linear (not necessarily continuous) map by sending any element of a Hamel basis of $A$ to any its preimage and extending it linearly to $A$.  Then $\phi$ can be defined as the $n$-th coordinate of this lift.
\end{remark}

Clearly every quasidiagonal trace is amenable and every MF-trace is hyperlinear. Whether the converse implications hold is an open problem.

% Daje esli b u nas ne bylo unitalnoy versii glavnoy teoremy, teorema pro sledy vse ravno poluchalas by, perehodya ot phi_n k phi_n^+  iz C(A^+) v multiplicatory M_n = M_n.

%One of biggest open problems is whether every amenable trace is quasidiagonal. TWW proved.........

\medskip

Throughout this paper we use notation $A^+$  for the minimal unitization of $A$.

\section{Contractive continuous sections}

For sets $B, D$ and a surjection $\pi: B \to D$, a {\it section} is a map $s: D\to B$ such that $s(d)\in \pi^{-1}(d)$, for any $d\in D$.

\medskip

It was proved in \cite{Bartle} that a continuous section always exists. Loring proved  that there exists a continuous section of norm arbitrary close to 1 (\cite{LoringSection}, \cite[Th. 3.2]{Manuilov}). Here we prove that there exists a continuous section of norm exactly 1. Results of this section will be used in the next ones.

\medskip

The following two lemmas are very well known. The second one is often used in lifting problems (see e.g. \cite{LoringBook}).

\begin{lemma}\label{NCinterval} Let $i, a, b\in A$  and  $0\le i \le 1$. Then $$\|(1-i)^{1/2}a(1-i)^{1/2} + i^{1/2}bi^{1/2}\| \le \max\{\|a\|, \|b\|\}. $$
\end{lemma}

\begin{lemma}\label{folklore}
Let $\pi: B\to B/I$ be a surjective $\ast$-homomorphism.
For any approximate unit  $\{u_{\lambda}\}$ in $I$ and any $x\in B$, one has $\limsup \|x(1- u_{\lambda})\| = \|\pi(x)\|$.
\end{lemma}

\medskip

Below we will use notation $B_{\delta}(x)$ for the open ball of radius $\delta$ centered at $x$, $S^1(A)$ for the unit sphere of a C*-algebra $A$,  and $A_{+}$ for the set of positive elements of $A$.

\medskip

\begin{theorem}\label{section} Any surjective $\ast$-homomorphism $\pi: B \to D$ has a continuous section $s: D \to B$ which is
\begin{enumerate}

\item  contractive,

\item  unital, if $\pi$ if unital,

\item selfadjoint (that is $s(a^*) = s(a)^*$ for any $a\in D$),

\item positive (that is $s(a)\ge 0$ for any $a\ge 0$).
\end{enumerate}

\noindent If $G$ is  a compact group, $B, D$ are $G$-algebras and $\pi$ is equivariant, then $s$  can be chosen equivariant.
\end{theorem}
\begin{proof}
First we will define  a section $s$ on the elements of $S^1(D)$. For that we define  a set-valued function $\Gamma: S^1(D)\to B$ by
$$\Gamma(x) = \begin{cases} 1, \; x=1 \\ \pi^{-1}(x) \bigcap S^1(B) \bigcap B_{+}, \; x\ge 0 \; \text{and} \; x\neq 1 \\  \pi^{-1}(x) \bigcap S^1(B), \; \text{otherwise}
\end{cases}.$$ Recall that Michael Selection Theorem would guarantee the existence of a continuous selection for $\Gamma$ if the following conditions held:

1) $S^1(D)$ is a paracompact space;

2) $B$  is a Banach space;

3) $\Gamma$  is lower semicontinuous;

4) for all $x\in S^1(D)$, the set $\Gamma(x)$ is nonempty, convex and closed.

 \noindent We observe that 1) holds because by Stone Theorem any metric space is paracompact,  and 2) is clear.
 Since for any $x\in S^1(D)$  one has $$\pi^{-1}(x) \bigcap S^1(B) = \pi^{-1}(x) \bigcap B_1(0),$$ the set $\Gamma(x)$ is convex.  As  is well known, a contraction lifts to a contraction, hence $\Gamma(x)$ is nonempty and thus 4) is also clear. We need to check 3).

 By definition of lower semicontinuity, $\Gamma$ is lower semicontinuous  at the point  $a$ if for any open set $V$ intersecting $\Gamma(a)$ there exists a neighbourhood $U$ of $a$ such that
$\Gamma(x)$ intersects  $V$ for all $x\in U$. So let us fix $a$ and $V$ and let $b_0 \in V\bigcap \Gamma(a)$.  There is a $\delta > 0$ such that $V\supset \mathcal B_{\delta}(b_0)$. There is a  $\delta' < \delta/2$ such that if $a\neq 1$, then $1\notin B_{\delta'}(a)$, and if $a\notin D_{+}$, then
$B_{\delta'}(a) \bigcap D_{+} = \emptyset$ (the latter is possible since $D_{+}$ is closed).
%$$\|z_1-z_2\|\le \delta' \; \Rightarrow  \||z_1|-|z_2|\|\le \delta/2$$ for any contractions $z_1, z_2$.

Let $U = \mathcal B_{\delta'}(a)\bigcap S^1(D)$. Let $\{i_{\lambda}\}$ be an approximate unit for $ker \;\pi$ quasicentral for $B$. Let  $x\in U$, $x\neq a$. Let $y\in S^1(B)$ be a preimage of $x$ which is positive if $x\ge 0$, and arbitrary otherwise.  By quasicentrality and Lemma \ref{folklore} there is a   sufficiently large   $\lambda$ so that
$$\|[(1-i_{\lambda})^{1/2}, y]\|< \delta'/3, \; \|[i_{\lambda}^{1/2}, b_0]\|< \delta'/3, \; \|(y-b_0)(1-i_{\lambda})\|< \|x-a\| + \delta'/3.$$
Let
$$\tilde y =  (1-i_{\lambda})^{1/2}y(1-i_{\lambda})^{1/2} + i_{\lambda}^{1/2}b_0i_{\lambda}^{1/2}.$$  Then $\tilde y$ is a preimage of $x$ and is positive if $x$ is positive. By Lemma \ref{NCinterval} $\|\tilde y\| \le 1$. Then, since $\tilde y$ is a preimage of $x$,  $\|\tilde y\| = 1$. Thus  $\tilde y \in \Gamma(x)$. We have
$$\|\tilde y - b_0\| = \|(y-b_0)(1-i_{\lambda})\| + 2\delta'/3 < \|x - a\| + \delta' \le \delta.$$ Thus $\tilde y \in \mathcal B_{\delta}(b_0)\subset V$. Therefore $\Gamma(x)$ intersects $V$ and the condition 3) is checked.

By  Michael Selection Theorem there is a continuous map $\tilde s:  S^1(D)\to B$ such that $\tilde s(x) \in \Gamma(x)$. Now we extend $\tilde s$ to $D$ by
$$\tilde s(x) = \begin{cases} \tilde s \left(\frac{x}{\|x\|}\right) \|x\|, \; x\neq 0\\0, \; x=0.\end{cases}$$
We obtained a contractive positive continuous section. Let $$s(x) = \frac{s(x)+s(x^*)^*}{2}, $$ $x\in D$. Then $s$ is a contractive positive continuous section which is also selfadjoint.

To obtain the equivariant version, one replaces  Michael Selection Theorem by its equivariant version \cite[Th. 3.4]{EquivariantSelection}. 
\end{proof}

\begin{remark} Applying Theorem \ref{section} to the standard action of circle  we can make our section homogeneous  (that is $s(\lambda a) = \lambda s(a)$).
\end{remark}

%It is well-known that each  asymptotic homomorphism is equivalent to an equicontinuous one. One does it using Using Bartle-Graves selection theorem (see e.g. \cite{Dadarlat2}).  It follows from Theorem \ref{section} that is addition one can make this equicontinuous asymptotic homomorphism to be contractive.

Here  are two immediate corollaries.

\begin{corollary}\label{asHomRquivContr} Every (discrete) asymptotic homomorphism is equivalent to an equicontinuous contractive self-adjoint positive (discrete) asymptotic homomorphism.
\end{corollary}
\begin{proof} An asymptotic homomorphism  $f_{\lambda}: A \to B$, $\lambda\in \Lambda$,   defines an actual $\ast$-homomorphism $f: A \to C_b(\Lambda, B)/C_0(\Lambda, B)$. By Theorem \ref{section} there is a contractive continuous self-adjoint positive section $s: C_b(\Lambda, B)/C_0(\Lambda, B) \to C_b(\Lambda, B)$. Then $ev_{\lambda}\circ s\circ f$, $\lambda\in \Lambda$, is an asymptotic homomorphism that  has all the required properties.
\end{proof}

\begin{corollary} In the definition of a hyperlinear trace (MF-trace, respectively) the maps $\phi_n$ can be chosen to be continuous, contractive, self-adjoint and positive.
\end{corollary}
\begin{proof} Proof is similar to the proof of  Corollary \ref{asHomRquivContr}, using the  quotient C*-algebra $\prod M_{k_n}/\oplus_2 M_{k_n}$ ($\prod M_{k_n}/\oplus M_{k_n}$, respectively).
\end{proof}

\medskip

The next theorem states that given a continuous section on a closed subset of a C*-algebra  we always can extend it to the whole C*-algebra.
It covers  \cite[Lemma 2.6.]{FGT} as a particular case.

\begin{theorem}\label{ExtendingSection} Let $\pi: B \to D$ be a surjective $\ast$-homomorphism, $E$ a closed subset of $D$, $s_0: E \to \pi^{-1}(E)$ a continuous (unital, contractive, self-adjoint, positive) section. Then there is a continuous (unital, contractive, self-adjoint, positive) section $s: D \to B$ such that $s|_E = s_0$.
\end{theorem}
\begin{proof} We will prove the existence of a continuous section extending $s_0$. The additional properties of unitality etc. one obtains similarly, by modifying the construction as in the proof of Theorem \ref{section}.

We define  a set-valued function $\Gamma: D\to B$ by
$$\Gamma(x) = \begin{cases} s_0(x), \; x\in E \\ \pi^{-1}(x), \; x\notin E.\end{cases}$$ We will show that $\Gamma$ is lower semicontinuous.  Let $a\in D$ and let $V$ be an open set intersecting $\Gamma(a)$. Let $b_0\in \Gamma(a)$. Since $V$ is open, there is a $\delta>0$ such that $B_{\delta}(b_0)\subset V$.

If $a\in E$, then, since $s_0$ is continuous, there is a $\delta'>0$ such that for any $x\in B_{\delta'}(a)\bigcap E$,
\begin{equation}\label{ExtendingSection1}s_0(x)\in B_{\delta}(s_0(a)) \subset V.\end{equation}

If $a\notin E$, then, since $E$ is closed, there is a $\delta'>0$ such that \begin{equation}\label{ExtendingSection2}B_{\delta'}(a)\bigcap E = \emptyset.\end{equation}

Either way let $\delta'' = \min\{\delta/2, \delta'\}$ and $U = B_{\delta''}(a).$ Let $x\in U$.

If $x\in E$, then by (\ref{ExtendingSection1}) $\Gamma(x) = s_0(x) \in V$.

If $x\notin E$, let  $\{i_{\lambda}\}$ be an approximate unit in $ker \;\pi$ quasicentral for $B$ and $y$ be any preimage of $x$. Let
$$\tilde y =  (1-i_{\lambda})^{1/2}y(1-i_{\lambda})^{1/2} + i_{\lambda}^{1/2}b_0i_{\lambda}^{1/2},$$ where $\lambda$ is sufficiently large to ensure that (by the same computation as in the proof of Th.\ref{section}) $$ \|\tilde y - b_0\|< \delta.$$ Thus $\tilde y\in B_{\delta}(b_0) \subset V$. On the other hand,
$\tilde y \in \pi^{-1}(x)$, and the latter coincides with $\Gamma(x)$ by (\ref{ExtendingSection2}).

Either way, we obtained that $\Gamma(x)$ intersects $V$, and hence $\Gamma$ is lower semicontinuous. Since the other assumptions of Michael Selection Theorem hold, we conclude that there is a continuous section $s$ which extends $s_0$ by construction.
\end{proof}

\section{Cones and applications}

\subsection{Cones and lifting to asymptotic homomorphisms}

\begin{lemma}\label{AsHom1}  Let $p (x_1, \ldots, x_N) $ be a homogeneous NC $\ast$-polynomial (more generally, $p$ can be of more variables and homogeneous only in $x_1, \ldots, x_N$). Let $I \lhd B$, $\{i_{\lambda}\}$ a quasicientral approximate unit of $I$ relative to $B$, and $\pi: B \to B/I$  the canonical surjection.
 Suppose $\pi(p(b_1, \ldots, b_N)) = 0$. Then $$\lim \|p(b_1(1-i_{\lambda}), \ldots, b_N(1-i_{\lambda}))\| =0.$$
\end{lemma}
\begin{proof} Let $d$ be the degree of homogeneuity of $p$. By quasicentrality of  $\{i_{\lambda}\}$ we have $$ \|p(b_1(1-i_{\lambda}), \ldots, b_N(1-i_{\lambda}))\| \approx \|p(b_1, \ldots, b_N) (1- i_{\lambda})^d\|.$$ Since $\{ 1- (1- i_{\lambda})^d\}$ is itself a (quasicentral) approximate unit, by Lemma \ref{folklore}  we obtain
$$\limsup \|p(b_1(1-i_{\lambda}), \ldots, b_N(1-i_{\lambda}))\| = \| \pi(p(b_1, \ldots, b_N)) \| =0.$$
\end{proof}

%\begin{lemma}\label{AsHom2} Let   $p (x_1, \ldots, x_N) $ be a NC $\ast$-polynomial which is homogeneous of degree $d$. Then there is a constant $c$ such that for any $a_1, \ldots, a_N, y$ with $\|a_i\|\le 1$  and $\|y\|\le 1$ $$\|p(y^{\frac{1}{2}}a_1y^{\frac{1}{2}}, \ldots, y^{\frac{1}{2}}a_Ny^{\frac{1}{2}}) - p(a_1, \ldots, a_N)y^d\| \le c \max_{1\le i\le N} \|[a_i,  y]\|^{\frac{1}{2}}.$$ \end{lemma} \begin{proof} One checks by induction that there is constant $c$ such that  $$\|p(y^{\frac{1}{2}}a_1y^{\frac{1}{2}}, \ldots, y^{\frac{1}{2}}a_Ny^{\frac{1}{2}}) - p(a_1, \ldots, a_N)y^d\| \le c \max_{1\le i\le N} \|[a_i,  y^{\frac{1}{2}}]\|.$$ The inequality $\|[x, y^{\frac{1}{2}}]\| \le \sqrt{2} \|[x, y]\|^{\frac{1}{2}}$ (\cite[Lemma 2.1]{OlsenPedersen}) finishes the proof. \end{proof}

In what follows we will consider universal C*-algebras generated by generators $x_1, x_2, \ldots$ and relations of the form $R_j(x_1, x_2, \ldots)=0$ or $\|R_j(x_1, x_2, \ldots)\|\le c$, $j\in \mathbb N$, where each $R_j$ is a noncommutative $\ast$-polynomial depending on finitely many variables.

We will use notation $\mathbf x = \{x_1, x_2, \ldots\}$ for the set of generators (also $\mathbf X$ for a tuple $\{X_1, X_2, \ldots\}$ in some C*-algebra $A$) and $\mathbf R = \{R_1, R_2, \ldots\}$ for the set of relations. We will write $C^*\langle \mathbf x\;|\;\mathbf R\rangle $ for the corresponding universal C*-algebra. We will write $\mathbf R(\mathbf X) =0$ meaning that $X_1, X_2, \ldots$ satisfy all the relations.

Below we let $C_b([1, \infty), B)$ denote the C*-algebra of all bounded continuous $B$-valued functions on $[1, \infty)$ and let $C_0([1, \infty), B)$ be the ideal of all functions vanishing at infinity.

\begin{lemma}\label{ExtendingFromGenerators} (i) Let $\phi: C^*\langle \mathbf x\;|\;\mathbf R\rangle \to B/I$ be a $\ast$-homomorphism and let $\mathbf X \in C_b([1, \infty), B)$  be such that
\begin{equation}\label{ExtendingFromGenerators1} \lim_{t\to\infty} \mathbf R(\mathbf X)(t) = 0\end{equation}
and
\begin{equation}\label{ExtendingFromGenerators1'} \pi(\mathbf X(t)) = \phi(\mathbf x), \end{equation}
for any $t\in [1,\infty)$.
Then there exists a contractive positive asymptotic homomorphism $f_t: C^*\langle \mathbf x\;|\;\mathbf R\rangle\to B$ such that $\pi\circ f_t = \phi$, for any $t\in [1,\infty)$, and $\lim_{t\to \infty} \|f_t(\mathbf x) - \mathbf X(t)\| = 0$.

\medskip

(ii) Let $p_1, p_2, \ldots$ be noncommutative $\ast$-polynomials.
 Let $B_0\subset B$ be a C*-subalgebra and suppose
$p_k(\mathbf X)\in C_b([0, \infty), B_0)$, for each $k\in \mathbb N$. Then the asymptotic homomorphism $f_t$ in (i) can be chosen with the additional property that $f_t(C^*(p_1(\mathbf x), p_2(\mathbf x), \ldots))\subset B_0$, for each $t$.
\end{lemma}
\begin{proof} (i)   From (\ref{ExtendingFromGenerators1'}) we obtain
$$\pi(\mathbf R(\mathbf x)(t)) = \mathbf R(\phi(\mathbf x)) =0,$$ for any $t$. Therefore $\mathbf R(\mathbf X)(t) \in I$ for any $t$ which together with (\ref{ExtendingFromGenerators1}) gives $\mathbf R(\mathbf X)\in C_0([1, \infty), I)$. Hence there is a $\ast$-homomorphism $$\Phi: C^*\langle \mathbf x\;|\;\mathbf R\rangle \to C_b([1, \infty), B)/C_0([1, \infty), I)$$ such that $\Phi(\mathbf x) = \mathbf X + C_0([1, \infty), I)$. Let
$s: C_b([1, \infty), B)/C_0([1, \infty), I) \to C_b([1, \infty), B)$ be a continuous contractive positive section guaranteed by Theorem \ref{section}. Define $f_t: C^*\langle \mathbf x\;|\;\mathbf R\rangle \to B$ as $$f_t = ev_t\circ s \circ \Phi.$$

Let $p$ be any noncommutative $\ast$-polynomial. Since $s$ is a section,
$$ s\circ \Phi(p(\mathbf x)) + C_0([1, \infty), I)   = \Phi(p(\mathbf x)) = p(\mathbf  X +  C_0([1, \infty), I)) = p(\mathbf X) + C_0([1, \infty), I).$$ Thus
$$s\circ \Phi(p(\mathbf x)) - p(\mathbf X) \in C_0([1, \infty), I).$$
Hence
\begin{equation}\label{ExtendingFromGenerators*} f_t(p(\mathbf x)) - p(\mathbf X)(t) \in I, \;\text{for any}\; t \end{equation}
and
\begin{equation}\label{ExtendingFromGenerators**} \lim_{t\to \infty} \|f_t(p(\mathbf x)) - p(\mathbf X)(t)\| = 0. \end{equation}
From (\ref{ExtendingFromGenerators*}) and (\ref{ExtendingFromGenerators1'}) we obtain
$$\pi(f_t(p(\mathbf x))) = \pi(p(\mathbf X)(t)) = p(\phi(\mathbf x)) = \phi(p(\mathbf x)),$$ for any $t$.  Thus for a dense subset of $C^*\langle \mathbf x\;|\;\mathbf R\rangle$,
$$\pi\circ f_t(a) = \phi(a),$$ for any $t$. Since for each $t$, $f_t$ is continuous by the construction, it follows that $\pi\circ f_t = \phi$, for any $t$.
Let us show that $f_t$ is an asymptotic homomorphism.
From (\ref{ExtendingFromGenerators**}) one deduces that for any NC polynomials $p, q$
\begin{equation}\label{ExtendingFromenerators0} \lim_{t\to\infty} \|f_t(p(\mathbf x)q(\mathbf x)) - f_t(p(\mathbf x))f_t(q(\mathbf x))\| =0.\end{equation} Let now $a, b \in C^*\langle \mathbf x\;|\;\mathbf R\rangle$, $\epsilon > 0$. Since NC polynomials of $\mathbf x$ are dense in
$C^*\langle \mathbf x\;|\;\mathbf R\rangle$ and $s\circ \Phi$ is continuous, there exist NC polynomials $p, q$ such that
$$\|s\circ \Phi(a) - s\circ \Phi(p(\mathbf x))\|<\epsilon, \; \|s\circ \Phi(b) - s\circ \Phi(q(\mathbf x))\|<\epsilon, \;
\|s\circ \Phi(ab) - s\circ \Phi(p(\mathbf x)q(\mathbf x))\|<\epsilon.$$ Therefore
\begin{equation}\label{ExtendingFromenerators00}\|f_t(a) - f_t(p(\mathbf x))\|< \epsilon, \;\|f_t(b) - f_t(q(\mathbf x))\|< \epsilon, \;\|f_t(ab) - f_t(p(\mathbf x)q(\mathbf x))\|< \epsilon,\end{equation}
for all $t$.
Using (\ref{ExtendingFromenerators0}), (\ref{ExtendingFromenerators00}) and contractibility of $f_t$ one easily checks that
$$\|f_t(ab) - f_t(a)f_t(b)\|< 4 \epsilon.$$ Thus $f_t$ is asymptotically multiplicative. Similarly one checks asymptotic linearity and self-adjointness. Also, for each $a\in A$, the map $t\mapsto f_t(a)$ is continuous by the construction.  Therefore  $f_t$ is an asymptotic homomorphism. It is contractive and positive by the construction, and $\lim_{t\to \infty} \|f_t(\mathbf x) - \mathbf X(t)\| = 0$  by (\ref{ExtendingFromGenerators**}).

\medskip

(ii)  We observe that since $p_k(\mathbf X)\in C_b([0, \infty), B_0)$, for each $k\in \mathbb N$,  the homomorphism $\Phi$ constructed in the proof of (i) sends the C*-subalgebra generated by $p_1(\mathbf x), p_2(\mathbf x), \ldots$ to the C*-subalgebra $C_b([0, \infty), B_0)/C_0([0, \infty), I\bigcap B_0) $ of $C_b([1, \infty), B)/C_0([1, \infty), I) $.  Let $s_0: C_b([1, \infty), B_0)/C_0([1, \infty), I\bigcap B_0) \to C_b([1, \infty), B_0)$ be a continuous contractive positive section guaranteed by Theorem \ref{section}. By Theorem \ref{ExtendingSection} we extend it to a continuous contractive positive section $s: C_b([1, \infty), B)/C_0([1, \infty), I) \to C_b([1, \infty), B)$. The rest of the proof of (i) goes without any change.
\end{proof}

\begin{remark}\label{DiscreteVersionOfLemma} In Lemma \ref{ExtendingFromGenerators} the parameter $t$ can be arbitrary. In particular Lemma \ref{ExtendingFromGenerators} applies also to obtain discrete asymptotic homomorphisms. In this case in the proof we replace $C_b([1, \infty)$ and $C_b([1, \infty)/C_0([1, \infty) $ by $\prod_{\lambda}B$ and $\prod_{\lambda}B /\oplus_{\Lambda} B$ respectively. \end{remark}

The following lemma is straightforward.

\begin{lemma}\label{continuousqau} Let $u_n, n\in \mathbb N,$ be a quasicentral approximate unit for $I\triangleleft B$. Let $u_t = (n+1-t)u_n + (t-n)u_n$, for $n\le t\le n+1$, $n\in \mathbb N$. Then the continuous path $u_t$, $t\in [1, \infty)$, is also a quasicentral approximate unit for $I\triangleleft B$.
\end{lemma}

\begin{theorem}\label{cones} Let $A$ be a  separable C*-algebra. Then any $\ast$-homomorphism from $CA$ to any quotient $B/I$ lifts to a contractive positive asymptotic homomorphism.
\end{theorem}
\begin{proof}  1) First we consider the case of unital $A$. By (\ref{ConePresentationTAMS}), $CA$ has presentation
$$CA = \langle h,  x_1, x_2, \ldots \;|\;  0\le h\le 1, -Ch \le x_i\le Ch, \; p_j (h, x_1, x_2, \ldots) = 0, \; i, j\in \mathbb N\rangle,$$ where each NC $\ast$-polynomial $p_j$ is homogeneous and depends on finitely many of $x_i$'s and their adjoints and possibly on $h$. For simplicity of notation we will assume we have only one polynomial $p$ and one variable $x$. The proof of general case is the same, one just handles the first $n$ polynomials at the $n$-th step of induction.  So WLOG
$$CA = \langle x, h \;|\; 0\le h\le 1, -Ch \le x \le Ch, \; p(h, x) = 0\rangle.$$ Let $d$ be the degree of homogeneuity of $p$.
 %and $c$ the constant from Lemma \ref{AsHom2}.

Let $f: A \to B/I$ be a $\ast$-homomorphism.
Lift $f(h)$ to $0\le H\le 1$. By Davidson's  two-sided order lifting theorem \cite[Lemma 2.1]{DavidsonLiftingPositive} we can lift $f(x)$ to $b\in B$ with  $-CH \le b\le CH$.

Let $\{i_{t}\}_{t\in [1, \infty)}$ be a quasicentral approximate unit for $I\triangleleft B$ as constructed in Lemma \ref{continuousqau}.
Let $$b_t= (1-i_t)^{1/2} b (1-i_t)^{1/2}, \; H_t =(1-i_t)^{1/2}H(1-i_t)^{1/2},$$
$t\in [1, \infty)$. Then $b_t$,$H_t$,  $t\in [1, \infty)$, are continuous paths and  $\pi(b_t)= x$ and $\pi(H_t)=h$, for any $t\in [1, \infty)$.
We also have $0\le H_t\le 1, -CH_t \le b_t\le CH_t$, and, by  Lemma \ref{AsHom1},  $ \lim_{t\to \infty}  \|p(H_t, b_t)\| = 0.$
Lemma \ref{ExtendingFromGenerators}  completes the proof of this case.

%Let $C_b([1, \infty), B)$ denote the C*-algebra of all bounded continuous $B$-valued functions on $[1, \infty)$ and let $C_0([1, \infty), B)$ be the ideal of all functions vanishing at infinity. Then $p(H_t, b_t) \in C_0([1, \infty), B)$ and therefore we can define a $\ast$-homomorphism  $\Phi: CA \to C_b([1, \infty), B)/C_0([1, \infty), B)$ by $$\Phi(h)  = H_t + C_0([1, \infty), B), \;\; \Phi(x) = b_t + C_0([1, \infty), B).$$  By Theorem \ref{section} there is a contractive continuous section  $s: C_b([1, \infty), B)/C_0([1, \infty), B) \to  C_b([1, \infty), B)$. Then $f_t:= ev_t\circ  s\circ \Phi: CA \to B$ is a contractive asymptotic homomorphism and for each $t$ we have $\pi\circ f_t(x) = \pi(b_t) = f(x)$ and $\pi\circ f_t(h) = \pi(H_t) = f(h)$, so each $f_t$ is a lift of $f$.

\medskip

2) Now assume $A$ is non-unital. In this case we need to slightly modify the arguments from 1). First, we can assume $f$ is surjective (by letting $B_0:= q^{-1}(f(CA))$ and considering $f$ as a homomorphism to $B_0/I$).
Since $CA$ is an essential ideal in $C(A^+)$, we have $CA \subset C(A^+)\subset M(CA)$. By the NC Tietze Extension Theorem $f$ extends to a homomorphism $ f': M(CA)\to M(B/I)$. Now we apply the arguments from 1) to
$\tilde f:=  f'|_{C(A^+)}: C(A^+)\to M(B/I).$
$$ \begin{tikzcd}  & & B
   \arrow[r,symbol=\subset]
 \arrow{d} & M(B) \arrow{d} \\ C(A^+)
 \arrow[bend left=30, dashed, swap]{urrr}
  \arrow[r,symbol=\supset]
  & CA \arrow{r}{f} & B/I
   \arrow[r,symbol=\subset]
  & M(B/I)
 \end{tikzcd}$$
We will use that

a)  as noted in Preliminaries, all the   elements $x_i$ and $x_ih^k$, $k\in \mathbb N$,  belong to $CA$ and generate it,

b) since $B$ is an ideal in $M(B)$, $b$ can be chosen to be in $B$ by the proof of \cite[Lemma 2.1]{DavidsonLiftingPositive}.

\noindent Then our construction gives $b_t\in B$, for each $t\in [1, \infty)$, and therefore $b_tH_t^k\in B$, for each $t\in [1, \infty)$, $k\in \mathbb N$. By Lemma \ref{ExtendingFromGenerators} (ii), there exists an asymptotic lift $f_t$ of $\tilde f$  such that $f_t|_{CA}$ lands in $B$.
\end{proof}

The following corollary shows that in Theorem \ref{cones} a cone can be replaced by any contractible C*-algebra.

\begin{corollary}  Let $A$ be a  separable contractible C*-algebra. Then any $\ast$-homomorphism from $A$ to any quotient $B/I$ lifts to a contractive positive asymptotic homomorphism.
\end{corollary}
\begin{proof}  Since $A$ is contractible, there is a $\ast$-homomorphism $\theta: A \to CA$ such that $ev_1\circ \theta = id_A$.  Let $f: A \to B/I$ be a $\ast$-homomorphism.
By Theorem \ref{cones}  $f\circ ev_1$ lifts to an asymptotic homomorphism $\psi_t: CD \to B$. Then $\psi_t\circ \theta$ is a asymptotic lift of $f\circ ev_1 \circ \theta = f$.
\end{proof}

\subsection{An application: Lifting to asymptotically cp maps}

%{\bf Obyasnit order 0 and as. order 0!!!!!!!!!}

In \cite{FGT} Gardella, Forough and Thomsen defined a continuous path $\theta_t:A \to B$ to be asymptotically positive, if each $\theta_t$  is self-adjoint and $\lim_{t\to \infty}\theta_t(a)^- = 0$. Asymptotic complete positivity they defined accordingly. They proved that any ccp map (order zero, respectively) $A\to B/I$ admits an asymptotically contractive asymptotically completely positive (and asymptotically order zero, respectively) lift, and obtained an equivariant version of it.

Here we obtain short proofs of their results with some improvements. We show that for an order zero map a lift can be chosen as a composition of an actual order zero map and an asymptotic homomorphism. For a ccp map a lift can be chosen as a composition of an order zero map and a corner of an asymptotic homomorphism. Moreover these lifts always can be made contractive and positive.

We start with a straightforward but convenient reformulation of asymptotic complete positivity.

\begin{proposition} TFAE:
\begin{itemize}
\item A continuous path $\theta_t:A \to B$ is asymptotically completely positive;

\item The map $\Theta: A \to C_b([1, \infty), B)/C_0([1, \infty), B)$ induced by $\theta_t$ by the formula $a\mapsto \theta_t(a) + C_0([1, \infty), B)$ is completely positive.
\end{itemize}
\end{proposition}

\begin{remark}\label{ascpNew} It follows from the proposition above that if $\phi_t$ is asymptotically cp and $b_t$ is a continuous path of positive elements, then $b_t\phi_tb_t$ is also asymptotically cp. In particular, a corner of an asymptotic homomorphism is asymptotically cp.
\end{remark}

\begin{corollary}\label{ShortProof}

(i) Any  ccp map lifts to a  contractive positive asymptotically cp map;

(ii)  Any ucp map lifts to a unital positive asymptotically cp map.

\end{corollary}
\begin{proof} (i):   Let $\phi: A \to B/I$ be ccp.   By Kasparov-Stinespring Dilation Theorem \cite{KasparovDilation} $\phi$ dilates to a $\ast$-homomorphism $\Phi: A \to M(B/I \otimes K)$. By Theorem \ref{cones}   $ev_1\circ \Phi: CA \to M(B/I\otimes K)$ lifts to an asymptotic homomorphism $\Psi_t: CA \to M(B\otimes K)$. Let $\delta: A \to CA$ be defined by $\delta(a) = t\otimes a$. Then $(\Psi_t)_{11}\circ \delta$ is a contractive positive asymptotically cp lift of $\phi$.

(ii): Let  $\phi: A \to B/I$ be ucp. Let $\psi_t$  be an asymptotically ccp lift of $\phi$. We adopt a construction from \cite[Lemma 13.1.2]{BO}.  For $t\in [1, \infty)$ define a continuous function $h_t$ on $[0, 1]$  by
$h_t(x) = \max\{\frac{1}{t}, x\}$ and fix a state $\omega$ on $A$. We define
$$\tilde\psi_t(a) = h_t(\psi_t(1))^{-1/2}\left( \psi_t(a) +\omega(a) (h_t(\psi_t(1)) - \psi_t(1))\right) h_t(\psi_t(1))^{-1/2}.$$ This is a unital
positive asymptotically cp lift of $\phi$.
\end{proof}

Recall that a completely positive map $\theta$ is {\it order zero} if for any positive $a, b$ with $ab=0$ one has $\theta(a)\theta(b)=0$. If this property holds asymptotically, then a continuous path is {\it asymptotically order zero}  \cite{FGT}.

 Now we obtain the order zero result of \cite{FGT} and additionally we make the lift contractive and positive.
 One can obtain it directly from Theorem \ref{cones} using duality between order zero maps on $A$ and homomorphisms from $CA$, which is done below. Alternatively it follows from Proposition \ref{equivariant} and Corollary \ref{ShortProof} which does not use the duality.

\begin{corollary}\label{corollaryOrderZero} Any contractive order zero map $\phi: A \to B/I$ lifts to a contractive positive asymptotically order zero  asymptotically cp map. \end{corollary}
\begin{proof} Let $\delta: A \to CA$ be given by $\delta(a) = t\otimes a$. By \cite{WZ} there is a $\ast$-homomorphism $\rho_{\phi}: CA \to B/I$ such that $\rho_{\phi}\circ\delta  = \phi$. By Theorem \ref{cones} we can lift it to a contractive  positive asymptotic homomorphism $\psi_t$.  Then $\psi_t\circ \delta$ is a contractive positive asymptotically order zero asymptotically cp map and a lift of $\phi$.
\end{proof}

\subsubsection{Equivariant version}
Besides lifting  cp maps to asymptotically cp, in \cite{FGT} they also consider the case when  $A, B, I$ are $G$-algebras, where $G$ is a second countable locally compact group,  and they show that a cp map always has an asymptotic lift which is asymptotically $G$-equivariant. Moreover they arrange that the asymptotic lift also asymptotically preserves the norm of products. Below we show how one can obtain this from Corollary \ref{ShortProof}.
The only ingredient needed is Kasparov's result about the existence of a quasicentral approximate unit $u_{n}$
%that satisfies  $u_{n+1}u_n = u_n$, $n\in \mathbb N$, and
which is $G$-invariant in the sense that $\lim_{n\to \infty} \max_{g\in K}\|gu_n - u_n \|= 0$ for any compact subset  $K\in G$ (\cite{Kasparov2}).

\begin{proposition}\label{equivariant} If a ccp map $\phi: A \to B/I$ lifts to a (contractive) (positive) asymptotically cp map $\psi_t: A \to B$, then it also lifts to a (contractive) (positive) asymptotically cp map $\bar\psi_t: A \to B$ with the additional properties that
$$\lim_{t\to \infty}\|\bar \psi_t(a)\bar \psi_t(b)\| = \|\phi(a)\phi(b)\|, $$ for any $a, b \in A$, and
$$\lim_{t\to \infty} \max_{g, h \in K}\|g\bar\psi_t(a)-\bar\psi_t(ha)\| = \|g\phi(a)-\phi(ha)\|,$$ for any $a\in A$ and any compact $K\subset G$.
In particular $\bar\psi_t$ is asymptotically equivariant when $\phi$ is equivariant, and $\bar\psi_t$ is asymptotically order zero  when $\phi$ is order zero.
\end{proposition}
\begin{proof}
%Let $t_n \to \infty$. Let $\{a_1, a_2, \ldots\}$ be a dense subset in $A$. Since $G$ is locally compact and second countable, there is a basis $U_1, U_2, \ldots$ such that $K_i = \overline U_i$ is compact, $i\in \mathbb N$. Then any compact $K\subset G$ is covered by finitely many $K_i$'s.  By the mentioned above Kasparov's result (\cite{Kasparov2}) and Lemma \ref{folklore} there is a quasicentral approximate unit $\{i_n\}$ in $I$ such that\begin{enumerate}

%\item $\max_{g\in K_k}\|g i_n - i_n \|\le 1/n,$ for $k\le n$,

%\item $\|[(1- i_n)^{1/2}, \psi_{t_m}(a_k)]\| \le 1/n,$ for $m\le n+1$, $k\le n$,

%\item $\|(1- i_n)^{m} \psi_{t_l}(a_k)\psi_{t_s}(a_j)\| - \|\phi(a_k)\phi(a_j)\| \le 1/n$, for $m\in \{1, 2\}$, $k, j\le n$, $l, s\in \{n, n+1\}$,

%\item $\max_{g, h \in K_k}\|(1-i_n)(\psi_{t_n}(ga)- h\psi_{t_n}(a))\| - \|\phi(ga) - h\phi(a)\| \le 1/n$, for $k\le n$,

%\item $i_{n+1}i_{n} = i_{n}i_{n+1} = i_n$, $n\in \mathbb N$

%\end{enumerate}
%(in fact, for us it would even be sufficient to replace the last condition  by the weaker condition that $\|[i_{n+1}, i_{n}]\|\le 1/n$).

%Let $\bar \psi_{t_n} = (1-i_n)^{1/2} \psi_{t_n} (1-i_n)^{1/2}$. For $n\le t\le n+1$ let $\delta = t-n$ and let
%$$\bar \psi_t = (1-\delta)\bar\psi_{t_n} + \delta \bar\psi_{t_{n+1}}.$$

Let $u_n$, $n\in \mathbb N$, be Kasparov's quasicentral approximate unit and let $u_t$ be as in Lemma \ref{continuousqau}. Then $\lim_{t\to \infty} \max_{g\in K}\|gu_t - u_t \|= 0$ for any compact subset  $K\in G$. Let $$\bar \psi_{t} = (1-u_t)^{1/2} \psi_t (1-u_t)^{1/2}.$$ By Remark \ref{ascpNew} $\psi_t$ is an asymptotically cp map. Using Lemma \ref{folklore} it is straightforward to check that $\psi_t$ satisfies all the other requirements.
\end{proof}

\begin{remark} Kasparov's result about approximate units assumes that $I$ is an ideal in a separable C*-algebra $B$. In the proposition above we do not assume that $B$ is separable, however since $A$ is separable, Kasparov's result is still applicable. Indeed let $H$ be a dense countable subgroup of $G$. Then the separable C*-subalgebra $D =C^*(H \phi(A)) \subset B/I$ is a $G$-algebra. Let $\{d_1, d_2, \ldots\}$ be a dense subset in $D$ and let $b_i\in B$ be any preimage of $d_i$, $i\in \mathbb N$. Then the separable C*-subalgebra $B_0 = C^*(H C^*(B_0)) \subset B$ is a $G$-algebra. Now Kasparov's result can be applied to the separable C*-algebra $B_0$ and its ideal $B_0\cap I$.
\end{remark}

\begin{remark} Theorems \cite[Th.3.2, Th. 3.4]{FGT} that state that if there is a completely positive lift, then one can find also a completely positive asymptotically $G$-invariant lift follow also from the proof of Proposition \ref{equivariant}.
\end{remark}

%{\bf One can look at the property that all asymptotically cp maps are close to cp (which is sort of WSP in category of cp maps) but it turns out to be nothing but LP.}

%\begin{proposition} TFAE: 1) Any asymptotic (or approximate) c.p. $A\to B$ is close to c.p.2) $A$ has the  LP.\end{proposition}
%{\bf Proof to be added.......}

\subsection{Lifting $\ast$-homomorphisms from cones - a completely positive version}

\begin{theorem}\label{LiftingDuality}   Let $A$  be a separable C*-algebra, $X$ a locally compact space, $\xi\in C_0(X)$ a strictly positive element, and let $\delta: A \to A\otimes C_0(X)$ be defined by $\delta(a)= a\otimes \xi$, $a\in A$. Let    $\rho: A\otimes C_0(X)\to B/I$  be a cp map.  Then $\rho$  lifts to a cp map if and only $\rho\circ\delta$ lifts to a cp map.
\end{theorem}
\begin{proof} One direction is clear. For the other direction,  we assume that $\rho\circ\delta$ lifts to a cp map. Let $F$ be a finite subset of the unit ball of $A\otimes C_0(X)$ and  $\epsilon >0$. By considering functions on $X$ as functions on its one-point compactification we can assume $X$ is compact. We refer to elements of $A\otimes C_0(X)$ as to $A$-valued functions on $X$.   There exists $\delta_1>0$ such that for any $\eta\in F$
\begin{equation}\label{LiftingDuality01}\|\eta(x) - \eta(x')\|< \epsilon/4, \;\text{whenever}\; |x-x'|< \delta_1.\end{equation}
There exists $\delta_2>0$ such that \begin{equation}\label{LiftingDuality2}\|\xi(x) - \xi(x')\|< \epsilon^2/8, \;\text{whenever}\; |x-x'|< \delta_2.\end{equation} Let
$$I_0 = \{x\in X\;|\; \|\eta(x)\|\le \epsilon/4, \eta\in F, \text{and}\; |\xi(x)|\le\epsilon/4.\}$$
Let $I_1, \ldots, I_N$ be  an open cover of $X\setminus I_0$ with $diam\; I_i < \min\{\delta_1, \delta_2\}$, $i=1, \ldots, N$.

%There exists $\Delta>0$ (one can take $\Delta = \delta_1$) such that \begin{equation}\label{LiftingDuality2}\|\eta|_{[0, \Delta]}\| < \epsilon/4,\end{equation} for any $\eta\in F$.  Let $\delta_2 = \min\{\frac{\epsilon\Delta}{2}, \delta_1\}.$ Let $I_1, \ldots, I_N$ be intervals inside $[\Delta, 1]$ that cover $[\Delta, 1]$ and have length not larger than $\delta_2$.

Then for any $t_1, t_2\in I_k$, $k = 1, \ldots,  N$, we have
\begin{equation}\label{fraction} \frac{\xi(t_1)}{\xi(t_2)} = 1+ \frac{\xi(t_1) - \xi(t_2)}{\xi(t_2)} \le 1+ \frac{\epsilon^2/8}{\epsilon/4} =1+ \epsilon/2.
\end{equation}
Choose a point $t_i\in I_i$, $i = 1, \ldots, N$. Define  $\alpha: A\otimes C_0(X) \to \mathbb C^N\otimes A $ by
$$\alpha(f) = \left(f(t_i)\right)_{i=1}^N, $$ $f\in A\otimes C_0(X)$. Let $\{u_i\}$ be a partition of unity for the cover of $X\setminus I_0$ by the open sets $I_i$, $i = 1, \ldots, N$. Define a cp map $\beta: \mathbb C^N\otimes A \to A\otimes C_0(X)$ by
$$\beta\left(\left(a_i\right)_{i=1}^N\right) = \sum_{i=1}^N \frac{a_i}{\xi(t_i)}\otimes \xi u_i.$$

{\it Claim 1:} For any $\eta\in F$, $$\|\beta\circ\alpha (\eta) - \eta\|< \epsilon.$$

{\it Proof of Claim 1:}  Using (\ref{LiftingDuality01}),(\ref{LiftingDuality2}), (\ref{fraction}), for any $\eta\in F$ we obtain
\begin{multline*}  \|\beta\circ\alpha (\eta) - \eta\| = \|\sum_{i=1}^N \frac{\eta(t_i)}{\xi(t_i)}\otimes \xi u_i - \sum_{i=1}^N \eta u_i - \eta|_{I_0}\| \\ \le \epsilon/4 + \|\sum_{i=1}^N \left(\frac{\eta(t_i)}{\xi(t_i)}\otimes \xi u_i - \eta u_i\right)\| \\
=\epsilon/4 + \sup_{t\in X\setminus I_0} \|\sum_{i=1}^N \left(\frac{\eta(t_i)}{\xi(t_i)}\xi(t)u_i(t) - \eta(t)u_i(t)\right)\|\\
=\epsilon/4 + \sup_{t\in X\setminus I_0} \|\sum_{i:\; t\in I_i} \left(\frac{\eta(t_i)}{\xi(t_i)}\xi(t)-\eta(t)\right)u_i(t)\|\\
\le \epsilon/4 + \sup_{t\in X\setminus I_0} \sum_{i:\; t\in I_i}\left( \|\eta(t_i)(\frac{\xi(t)}{\xi(t_i)}-1)\|+ \|\eta(t_i)-\eta(t)\|\right)u_i(t) \\
\le \epsilon/4 + \sup_{t\in X\setminus I_0} \sum_{i:\; t\in I_i}\frac{3\epsilon}{4} u_i(t) = \epsilon/4 + \sup_{t\in X\setminus I_0} \frac{3\epsilon}{4}
\sum_{i:\; t\in I_i} u_i(t) \le \epsilon/4 + 3\epsilon/4  = \epsilon.
\end{multline*}
Claim 1 is proved.

\medskip

{\it Claim 2:} $\rho\circ \beta$, and therefore $\rho\circ \beta\circ \alpha$,   lifts to a cp map.

{\it Proof of Claim 2:}  Let
$$\beta_i(a) := a\otimes \xi u_i = \left(1\otimes u_i^{1/2}\right)(a\otimes \xi) \left(1\otimes u_i^{1/2}\right) = \left(1\otimes u_i^{1/2}\right)\delta(a) \left(1\otimes u_i^{1/2}\right),$$ $a\in A$. Then
$$\rho \circ \beta_i(a)  = \rho\left(\left(1\otimes u_i^{1/2}\right)\delta(a) \left(1\otimes u_i^{1/2}\right)\right) = \rho\left(\left(1\otimes u_i^{1/2}\right)\right) \rho\circ\delta(a) \rho\left(\left(1\otimes u_i^{1/2}\right)\right), $$ $a\in A$.
Since $\rho\circ\delta$ lifts to a cp map and $\rho\left(\left(1\otimes u_i^{1/2}\right)\right)$ lifts to a positive element, $\rho\circ \beta_i$ lifts to a cp map. Let $pr_i: \mathbb C^N\otimes A\to A$ be the projection map on the $i$-th coordinate. Since $\beta = \sum_{i=1}^N \frac{1}{\xi(t_i)}\beta_i\circ pr_i$, $\rho\circ \beta$ also lifts to a cp map.  Claim 2 is proved.

\medskip
Let $\{a_1, a_2, \ldots\}$ be a dense subset of the unit ball of $A$.
Taking $\epsilon = 1/n$ and $F= \{a_1, \ldots, a_n\}$, by Claim 1 we obtain that $\rho$ is pointwise limit of maps $\rho\circ \beta\circ\alpha$ which lift to cp  maps by Claim 2. By Arveson's classical result, $\rho$ lifts to a cp map.
\end{proof}

\begin{lemma}\label{cpLiftFromUnitization} \begin{itemize}
 \item[(1)] Let $0\to J\to E\to A\to 0$ be a  short exact sequence, where $E$ is separable  and  $A$ has the  Lifting Property (LP). Let  $\psi: E\to B/I$ be a cp map. If $\psi|_J$ has a cp lift, then $\psi$ has a cp lift.

\item[(2)] If $\psi: C(A^+)\to B/I$ is a cp map and $\psi|_{CA}$ has  a cp lift, then $\psi$ has  a cp lift.
\end{itemize}
\end{lemma}
\begin{proof} 1) Since $A$ has the LP, the quotient map $p: E \to A$ has a cp split $s: A\to E$. Let $\{u_{\lambda}\}$ be a quasicentral approximate unit for $J$ relative to $E$. Using Lemma \ref{folklore} it is straightforward to check that for any $x\in E$ we have
$$x = \lim_{\lambda\to \infty} \left(u_{\lambda}^{1/2}xu_{\lambda}^{1/2} + (1-u_{\lambda})^{1/2}  s\circ p(x) (1-u_{\lambda})^{1/2}\right).$$
Therefore
$$\psi(x) = \lim_{\lambda\to \infty} \left(\psi\left(u_{\lambda}^{1/2}xu_{\lambda}^{1/2}\right) + \psi\left((1-u_{\lambda})^{1/2}  s\circ p(x) (1-u_{\lambda})^{1/2}\right)\right), $$ for any $x\in E$. Define $\bar\psi_{\lambda}: E \to B/I$ and $\tilde \psi_{\lambda}: E \to B/I$ by
$$\bar\psi_{\lambda}(x) = \psi\left(u_{\lambda}^{1/2}xu_{\lambda}^{1/2}\right),$$ $$ \tilde\psi_{\lambda}(x) = \psi\left((1-u_{\lambda})^{1/2}s\circ p(x)(1-u_{\lambda})^{1/2}\right).$$ Then $\psi$ is the pointwise limit of the cp maps $\bar\psi_{\lambda} + \tilde \psi_{\lambda}$. Let $\phi$ be a cp lift of $\psi|_J$. Then $\bar\psi_{\lambda}$ has a cp lift $x\mapsto \phi(u_{\lambda}^{1/2}xu_{\lambda}^{1/2})$. Since $\tilde\psi_{\lambda}$ factorizes through $A$ and $A$ has the LP,  $\tilde\psi_{\lambda}$  also has a cp lift. Therefore $\psi$ is a pointwise limit of liftable cp maps and hence is liftable itself.

2) Apply 1) to the sequence $0\to CA\to C(A^+)\to C\mathbb C\to 0$.
\end{proof}

\begin{theorem}\label{ccp}   Let $A$ be separable and let $f: CA \to B/I$ be a $\ast$-homomorphism that has a cpc lift $\phi: CA \to B$.
 Then $f$ lifts to a cpc  asymptotic homomorphism.
  \end{theorem}
  \begin{proof} If $A$ is non-unital, we extend $f$ to $\tilde f: C(A^+)\to M(B/I)$ in the same way as in part 2) of the proof of Theorem \ref{cones}.
  By Lemma \ref{cpLiftFromUnitization}, $\tilde f$ lifts to a cpc map $ \psi: C(A^+)\to M(B)$. If $A$ is unital, we let $\psi = \phi$.

  Let $\{i_{\lambda}\}_{\lambda\in \Lambda}$ be  a quasicentral approximate unit of $I$ relative to $B$. Here $\Lambda=\mathbb N$ is we are interested in discrete asymptotic homomorphisms and $\Lambda = [1, \infty)$ for the continuous version in which case we take a quasicentral approximate unit forming a continuous path as in Lemma \ref{continuousqau}.
Let $$\psi_{\lambda} = (1-i_{\lambda})^{1/2} \psi (1-i_{\lambda})^{1/2}: C(A^+)\to M(B).$$ For each $\lambda\in \Lambda$ and $x\in CA$ we have $$q(\psi_{\lambda}(x)) = f(x).$$
 Let  $ a_1, a_2 \ldots$ be a dense subset of $A$. Let $x_i = \delta(a_i)$, $h = \delta(1_{A^+})$.
We write $C(A^+)$ as the universal C*-algebra with generators $h$ and $x_i$, $i\in \mathbb N$, and homogeneous relations.   By Lemma \ref{AsHom1}  $\psi_{\lambda}(h)$ and $\psi_{\lambda}(x_i)$'s approximately satisfy the relations.  Let $$\pi: \prod M(B) \to \prod M(B) /\oplus I$$ ($\pi: C_b([1, \infty), M(B))\to C_b([1, \infty), M(B))/C_0([1, \infty), I)$ for the continuous version) be the canonical surjection. Then $\pi\left((\psi_{\lambda}(h))_{\lambda\in \Lambda}\right)$ and  $\pi\left((\psi_{\lambda}(x_i))_{\lambda\in \Lambda}\right)$, $i\in \mathbb N$, satisfy all the relations of $C(A^+)$ and therefore define a $\ast$-homomorphism
$$\Psi: C(A^+) \to \prod M(B)/\oplus I$$ ($\Psi: C(A^+) \to C_b([1, \infty), M(B))/C_0([1, \infty), I)$ for the continuous version). Let
$$\Phi = \Psi|_{CA}.$$ We notice that it lands in $\prod B/\oplus I$ (in $C_b([1, \infty), B)/C_0([1, \infty), I)$ for the continuous version). We have
$$\Phi(\delta(a_i)) = \Phi(x_i) = \pi\left((\psi_{\lambda}(x_i))_{\lambda\in \Lambda}\right) =  \pi\left((\psi_{\lambda}(\delta(a_i)))_{\lambda\in \Lambda}\right). $$ Since $\{a_i\}$ is dense in $A$, we conclude that
$$\Phi\circ \delta =\pi\circ \left(\psi_{\lambda}|_{CA}\circ \delta\right)_{\lambda\in \Lambda}.$$
Thus $\Phi\circ \delta$ has a cp lift $\left(\psi_{\lambda}|_{CA}\circ \delta\right)_{\lambda\in \Lambda}$. (We note that one cannot say that $\left(\psi_{\lambda}|_{CA}\right)_{\lambda\in \Lambda}$ is a cp lift of $\Phi$. It only lifts $\Phi$ on the linear span of the generators). Then, by Theorem \ref{LiftingDuality}, $\Phi$ has a cp lift
 $$\left(f_{\lambda}\right)_{\lambda\in\Lambda}: CA \to \prod B$$ ($\left(f_{\lambda}\right)_{\lambda\in\Lambda}: CA \to C_b([1, \infty, B)$ for the continuous version. In this case we also obtain that the map $\lambda\mapsto f_{\lambda}$ is continuous).

  Then for any $x, y\in CA$ we have
 $$\left(f_{\lambda}(xy) - f_{\lambda}(x)f_{\lambda}(y)\right)_{\lambda\in \Lambda}\in \oplus I  \; $$
  ($\in (C_0([1, \infty), I), \;\text{respectively}$) and similar for linear combinations and adjoint elements. Therefore
 $f_{\lambda}$, $\lambda\in \Lambda$, is an asymptotic homomorphism such that
  $q\circ f_{\lambda}$ is a $\ast$-homomorphism for each $\lambda\in \Lambda$.

  By the construction of $\Phi$ and $f_{\lambda}$, for each generator $x_i$ of $CA$ we have
  $$\pi\left((f_{\lambda}(x_i))_{\lambda\in\Lambda}\right) = \Phi(x_i) = \pi\left((\psi_{\lambda}(x_i))_{\lambda\in\Lambda}\right).$$
  Hence $$q\circ f_{\lambda}(x_i) = q\circ \phi_{\lambda}(x_i) = f(x_i),$$ for each $\lambda\in \Lambda$. Therefore
  $q\circ f_{\lambda} = f, $ for each $\lambda \in \Lambda$.
  \end{proof}

\subsubsection{Byproduct of some results of this section}

As a byproduct of Theorem \ref{LiftingDuality} we obtain the following statement that can be of independent interest, and in particular a short proof to  \cite[Lemma 4.1]{BCW}. Recall that by Winter-Zacharias duality between order zero maps from $A$  and $\ast$-homomorphisms from $CA$, for any order zero map $\phi:A \to B$ there is a $\ast$-homomorphism $\rho_{\phi} : CA \to  B$ such that $\rho_{\phi}\circ \delta = \phi$ (\cite{WZ}).

\begin{corollary}\label{OrderZeroVsHomomorphism} Let $A$  be a separable C*-algebra and $\phi:A \to B/I$  an order zero map. If $\phi $ lifts to a cp map, then $\rho_{\phi}$ lifts to a cp map.
\end{corollary}
\begin{proof} Since $\phi = \rho_{\phi}\circ \delta$, the statement follows from Theorem  \ref{LiftingDuality}.
\end{proof}

\begin{corollary}\label{lemmaBCW}(\cite[Lemma 4.1]{BCW}) Let A be a separable unital C*-algebra, and let B be a unital
C*-algebra. Suppose that $(\psi_n)_{n}$ is a sequence of c.p.c. maps from $A$ into
$B$ inducing an order zero map $\psi : A \to B_{\omega}$. Then there exists a sequence of
cpc maps $(\phi_n)_n$ from $CA$ into $B$ inducing a $\ast$-homomorphism
$\phi : CA\to B_{\omega}$ such that $\phi(t\otimes x) = \psi(x)$, for $x\in A$.
\end{corollary}
\begin{proof} Apply  either Theorem \ref{LiftingDuality} or Corollary \ref{OrderZeroVsHomomorphism} to $B/I: =B_{\omega}$ and $B:= \prod B$.
\end{proof}

While the Local Lifting Property (LLP) is known to be preserved under extensions (\cite[Cor. 2.6.]{Kirchberg}), this seems to be unknown for the Lifting Property (LP). We obtain it here as a consequence of  Lemma \ref{cpLiftFromUnitization}.

\begin{corollary}\label{LP} If $0\to J\to E\to A\to 0$ is a short exact sequence,  $E$ is separable, and $J$ and $A$ have the Lifting Property (LP), then $E$ has the LP.
\end{corollary}

\medskip

\subsection{An application: Traces on cones}

Brown, Carrion and White proved that each amenable trace on a cone is quasidiagonal \cite{BCW}. We start with giving a new proof to this  result.

\begin{theorem}\label{QDtraceCone} (Brown, Carrion, White \cite[Prop. 3.2]{BCW})  Let $A$ be  a C*-algebra. Then any amenable trace on $CA$ is quasidiagonal.
\end{theorem}
\begin{proof} Since quasidiagonality of a trace  is a local property and by Arveson's Extension Theorem cp maps to matrices extend from subalgebras, it is sufficient to assume that $A$ is separable. Let $\tau$ be an amenable trace on $CA$. Then there are cpc maps $\phi_n: CA \to M_{k_n}$ such that
$$\|\phi_n(ab) - \phi_n(a)\phi_n(b)\|_2\to 0, \;\; |\tau(a)- tr \phi_n(a)|\to 0 \;\;\text{as} \;\;n\to \infty,$$
for any $a, b \in CA$.  Let $\phi = (\phi_n)_{n\in \mathbb N}: CA \to \prod M_{k_n}$.  Let  $q: \prod M_{k_n}\to \prod M_{k_n}/\oplus_2 M_{k_n}$ be the canonical surjection. Then
$f= q\circ \phi$ is a $\ast$-homomorphism. By Theorem \ref{ccp} there is an asymptotic homomorphism  $\psi^{\lambda}$, $\lambda\in \Lambda$, consisting of cpc lifts of $f$.
\noindent Since each $\psi^{(\lambda)}$ is a lift of $f$, for each $\lambda$ we have
\begin{equation}\label{trace3} \lim_{n\to \infty} \|\phi_n(a)- \psi_n^{(\lambda)}(a)\|_2 =0, \end{equation}
for any $a\in CA$.  Now let $G$ be a finite subset of $CA$ and $\epsilon > 0$. Since the family $\psi^{(\lambda)}$, $\lambda\in \Lambda$, is an asymptotic homomorphism, there is $\lambda_0$ such that
\begin{equation}\label{trace4} \|\psi^{(\lambda_0)}(ab) - \psi^{(\lambda_0)}(a)\psi^{(\lambda_0)}(b)\|< \epsilon, \end{equation} for $a, b \in G$.  Therefore for each $n\in \mathbb N$
\begin{equation}\label{trace5} \|\psi_n^{(\lambda_0)}(ab) - \psi_n^{(\lambda_0)}(a)\psi_n^{(\lambda_0)}(b)\|< \epsilon. \end{equation}
By (\ref{trace3}) there is $n_0$ such that for any $n>n_0$
$$\|\phi_n(a) - \psi_n^{(\lambda_0)}(a)\|_2< \epsilon/2, $$
$a\in G$. There is $n_1\in \mathbb N$ such that for any $n>n_1$
$$|\tau(a)- tr \phi_n(a)| < \epsilon/2,$$ $a\in G$. Then for any $n> \max\{n_0, n_1\}$ we obtain
\begin{multline}\label{trace6} |\tau(a) - tr \psi_n^{(\lambda_0)}(a)|\le |\tau(a) - tr \phi_n(a)| + | tr \phi_n(a) - tr \psi_n^{(\lambda_0)}(a)|\\ \le \epsilon/2 + \|\phi_n(a) - \psi_n^{(\lambda_0)}(a)\|_2 <\epsilon.\end{multline}
By (\ref{trace5}) and (\ref{trace6}) $\tau$ is quasidiagonal.
\end{proof}

\medskip

We also obtain a similar result for hyperlinear traces which seems to be new.

\begin{theorem}\label{MFtraceCone} Let $A$ be  a C*-algebra. Then any hyperlinear trace on $CA$ is MF.
\end{theorem}
\begin{proof} Same proof, only using Theorem \ref{cones} instead of Theorem \ref{ccp}.
\end{proof}

\subsubsection{Traces on suspensions}
As is well known, any suspension C*-algebra $SA = A\otimes C_0(0, 1)$ is QD (and therefore MF), just because it is a subalgebra of the cone C*-algebra $CA$. However nothing seems to be known about traces on suspensions.

\begin{proposition}\label{hypTraceExtention} Let $I$ be an ideal of $A$. Then any hyperlinear trace on $I$ extends to a hyperlinear trace on $A$.
\end{proposition}
\begin{proof} Let $\tau$ be a hyperlinear trace on $I$. By Remark \ref{hypTraceReformulation}, $\tau = tr \; f$, for some $\ast$-homomorphism $f: A \to \prod M_{k_n}/\oplus_{2, \omega} M_{k_n}$. The ultraproduct C*-algebra $\prod M_{k_n}/\oplus_{2, \omega} M_{k_n}$ is a von Neumann subalgebra of $B(H)$, for some Hilbert space $H$. Then $f$, being a $\ast$-homomorphism to $B(H)$, extends to a $\ast$-homomorphism $\tilde f: A \to B(H)$ by the well-known formula $$\tilde f(a) = SOT-\lim_{\lambda\to \infty} f(i_{\lambda}a),$$
$a\in A$. Here $\{i_{\lambda}\}$ is an approximate unit for $I$. Since $\prod M_{k_n}/\oplus_{2, \omega} M_{k_n}$ is closed in the strong operator topology in $B(H)$, $\tilde f$ is a $\ast$-homomorphism from $A$ to $\prod M_{k_n}/\oplus_{2, \omega} M_{k_n}$. Again by Remark \ref{hypTraceReformulation}, $\tilde \tau:= tr \;\tilde f$ is a hyperlinear trace on $A$, and it extends $\tau$.
\end{proof}

\begin{corollary}\label{suspension} Let $I$ be an  ideal of $A$ and suppose that all hyperlinear traces on $A$ are MF. Then all hyperlinear traces on $I$ are MF. In particular, for any $A$, all hyperlinear traces on $SA$ are MF.
\end{corollary}
\begin{proof} Let $\tau$ be a hyperlinear trace on $I$. By Proposition \ref{hypTraceExtention}, $\tau$ extends to a hyperlinear trace $\tilde \tau$ on $A$. By the assumption,  $\tilde \tau$ is then MF, that is
$\tilde\tau(a) = \lim_{n\to \infty} tr \; \phi_n(a)$, for all $a\in A$, where $\phi_n:A \to M_{k_n}$ are approximately multiplicative, approximately linear and approximately self-adjoint maps. Then  $\tau(i) = \lim_{n\to \infty} tr \; \phi_n|_I(i)$, for all $i\in I$. So $\tau$ is MF.

Since $SA$ is an ideal in $CA$, the statement about $SA$ follows from Theorem \ref{MFtraceCone}.
\end{proof}

%\begin{remark} Similar arguments do not work to show that any amenable trace on $SA$ is quasidiagonal, because it is not clear whether an analogue of Proposition \ref{hypTraceExtention} holds for amenable traces. \end{remark}

%\section{New proofs of some quasidiagonality results}

\subsection{An application: A new proof of quasidiagonality of cones}

Recall that a C*-algebra $A$ is {\it quasidiagonal} (QD) if for each $\epsilon > 0$ and finite $\mathcal F \subset A$ there exist $n\in \mathbb N$ and a cpc map $\phi: A \to M_n$ such that
$$\|\phi(ab) - \phi(a)\phi(b)\|\le \epsilon,$$
$$\|\phi(a)\|\ge  \|a\|- \epsilon,$$
for any $a, b\in \mathcal F$.

\begin{proposition}\label{Voiculescu} (Voiculescu \cite{Voiculescu}) The cone over any C*-algebra is QD.
\end{proposition}
\begin{proof}   It is sufficient to prove the statement for a separable C*-algebra $A$.  Let $H$ be a Hilbert space and let $P_{n}$, $n \in \mathbb N$,  be an increasing sequence of projections of dimension $n$ that $\ast$-strongly converge to $1_{B(H)}$. We will identify $M_{n}$ with $P_{n} B(H) P_{n}$.
Let $\mathcal D \subset \prod_{n\in \mathbb N} M_{n}$ be the C*-algebra of all $\ast$-strongly convergent sequences of matrices and let $q: \mathcal D \to B(H)$ be the surjection that sends each sequence  to its $\ast$-strong limit.
   Let $\rho: CA \to B(H)$ be an embedding. It has a cpc lift $\left(P_n\rho P_n\right)_{n\in \mathbb N}: CA \to \mathcal D$.
By Theorem \ref{ccp} there is an asymptotically multiplicative family $\psi^{\lambda}: CA \to \mathcal D$ of cpc lifts of $\rho$. Let $\mathcal F \subset CA$ be finite and $\epsilon >0$. Since $\psi^{\lambda}$ is an asymptotically multiplicative family, there is $\lambda_0$  such that
\begin{equation} \|\psi^{\lambda_0}(ab) - \psi^{\lambda_0}(a)\psi^{\lambda_0}(b)\|\le \epsilon, \end{equation} for any
$a, b \in \mathcal F$. Hence
\begin{equation}\label{QD1} \|\psi_n^{\lambda_0}(ab) - \psi_n^{\lambda_0}(a)\psi_n^{\lambda_0}(b)\|\le \epsilon, \end{equation} for any
$n\in \mathbb N$, $a, b \in \mathcal F$.
Since $\psi^{\lambda_0}$ is a lift of $\rho$, it is isometric. Hence there is $N\in \mathbb N$ such that
\begin{equation}\label{QD2} \|\psi_N^{\lambda_0}(a)\|\ge \|a\|- \epsilon,\end{equation} for any $a\in \mathcal F$.
By (\ref{QD1}) and (\ref{QD2}) $CA$ is QD.
\end{proof}

It is not a better proof than the original proof  of this fact (\cite{Voiculescu}), but it is different and  based only on the algebraic structure of cones.
%In particular it shows that it is not necessary to use Voiculescu's theorem to prove this result.

\end{document}